\documentclass[11pt]{article}
\usepackage[numbers]{natbib}
\usepackage{amsbsy}
\usepackage{amsgen}
\usepackage[small,nohug,heads=vee]{diagrams}
\usepackage{amstext}
\usepackage{amsxtra}
\usepackage{amsmath}
\usepackage{amsopn}
\usepackage{amsthm}
\newtheorem{thm}{Theorem}[section]
\theoremstyle{definition}
\newtheorem{dfn}{Definition}[section]
\theoremstyle{plain}

\theoremstyle{remark}
\newtheorem{note}{Remark}[section]
\theoremstyle{plain}

\theoremstyle{plain}

\theoremstyle{plain}

\theoremstyle{conjecture}

\begin{document}\date{}
\title{Stability Conditions and Mirror Symmetry of K3 Surfaces in Attractor Backgrounds}
\author{Wenxuan Lu\footnote{wenxuanl@math.upenn.edu}} \maketitle
\begin{abstract}We study the space of stability conditions on $K3$ surfaces from the perspective of mirror symmetry. It is done in the so called attractor backgrounds (moduli)  which can be far  from the conventional large complex limits and are selected by the attractor mechanism for certain black holes. We find  certain highly non-generic behaviors of stability walls (a key notion in the study of wall crossings) in the space of stability conditions. They correspond  via mirror symmetry to some non-generic behaviors of special Lagrangians in an attractor background. The main results can be understood as a mirror correspondence in a synthesis of homological mirror conjecture  and SYZ mirror conjecture.     \end{abstract}

\newpage

\tableofcontents{\section{Introduction}}

 We study in this paper the interaction of the following three topics on $K3$ surfaces\begin{enumerate}\item Stability conditions on derived categories of coherent sheaves. \item Mirror Symmetry.\item Attractor mechanism in the study of certain black holes.\end{enumerate}

The study of stability conditions is a generalization of the study of stable bundles. The space of stability conditions on derived categories of coherent sheaves is a generalization of the space of complexified Kahler moduli. The study of counting stable objects is closely  related to the enumerative geometry.  If we move in the space of stability conditions we might encounter real codimension one walls (called stability walls in this paper) such that the counting problem of stable objects exhibits nontrivial jumps called wall crossings. Moreover stability conditions are actually defined  for general triangulated categories which include scenarios beyond complex geometry.\\

Stability conditions is a key notion that connects the two major programs of mirror symmetry: SYZ mirror conjecture of Strominger, Yau and Zaslow and  homological mirror conjecture of Kontsevich.

Roughly speaking the homological mirror conjecture predicts that for a mirror pair of Calabi-Yau varieties $M$ and $\check{M}$ there should be a duality exchanging the complex geometry and symplectic geometry of them. We expect a correspondence of the following form\footnote{This is meant to be a very short   description of the background. So we have deliberately ignored  many important details.}\begin{equation}DFuk(M)\longleftrightarrow \mathcal{D}(\check{M})\end{equation}Here $DFuk(M)$ is the derived Fukaya category of $M$. It is obtained by applying some categorical construction to the Fukaya category. Lagrangian submanifolds are objects of the Fukaya category. $\mathcal{D}(\check{M})$ is the derived category of coherent sheaves on $\check{M}$.

On the other hand SYZ mirror conjecture predicts that $M$ carries a special Lagrangian torus fibration structure and $\check{M}$ should be constructed as the dual special Lagrangian fibration. Passing from Lagrangians to special Lagrangians means that we need to impose conditions on the so called central charge $Z(\gamma,\Omega)$ (see section 2 for the definition). It turns out when we
define stability conditions there is a similar procedure involving central charges. We can speculate that there are   stability conditions on  $DFuk(M)$ satisfying essentially the same axioms in the definition of stability conditions on derived categories of coherent sheaves. Moreover special Lagrangians should be realized as stable objects.

 On the other side of the mirror duality (i.e. $\mathcal{D}(\check{M})$) we will also put stability conditions on $\mathcal{D}(\check{M})$ and then consider stable objects. Mirror symmetry should then exchange stable objects of the two sides. In fact in some cases one can check the correspondence of special Lagrangians in $M$ and stable bundle on $\check{M}$ which can be interpreted as a baby version of this correspondence. Since we expect wall crossings when we count stable objects there should also be a correspondence of wall crossings of the two sides.\\

So we have explained that we expect a much more ambitious and unified version of mirror symmetry of the following form
\begin{equation}DFuk(M)+ stability\ conditions\longleftrightarrow \mathcal{D}(\check{M})+ stability\ conditions\end{equation}
\begin{equation}stable\ objects\ in\ DFuk(M)\longleftrightarrow stable\ objects\ in\ \mathcal{D}(\check{M})\end{equation}
$$space\ of\ stability\ conditions\ of\ DFuk(M)\longleftrightarrow$$\begin{equation}\longleftrightarrow space\ of\ stability\ conditions\ of\ \mathcal{D}(\check{M})\end{equation}
\begin{equation}wall\ crossings\ for\ DFuk(M)\longleftrightarrow wall\ crossings\ for\ \mathcal{D}(\check{M})\end{equation}

As far as we know there is very few work in the literature about this version of mirror symmetry because it seems extremely difficult to study stability conditions on derived Fukaya categories. In particular it seems that there is virtually no results about correspondences of type (5).\\

One of the purposes of this paper is to formulate and prove a few theorems that can be considered as a  softened version of a correspondences of type (5) for certain $K3$ surfaces. To explain it we need to introduce the third player of this paper: attractor backgrounds (moduli) produced by the attractor mechanism  for black holes.\\

In the study of certain four dimensional black hole solutions obtained by reduction from a ten dimensional string theory  physicists have discovered an interesting fact called the attractor  mechanism. We assume that the black hole solution has spherical symmetry so that it makes sense to talk about the radial direction. The black hole solution of course depends on the six dimensional part of the spacetime. It turns out the six dimensional part over each value of the radial direction should be a Calabi-Yau threefold and these threefolds can have different moduli. Now the attractor mechanism claims that the underlying supersymmetry of the theory forces the moduli points over the horizon of the black hole to be special. They are called attractor backgrounds in this paper. In fact there is a flow called the attractor flow in the moduli space and the attractor backgrounds are stationary points.

Attractor backgrounds are very special in the moduli space. If the Calabi-Yau threefold is $K3\times T^{2}$ one can determine them explicitly. It is then interesting to ask what mirror symmetry can say about them. In this paper we study the mirror symmetry of $K3$ surfaces coming from the $K3\times T^{2}$ attractor backgrounds. This is  somewhat unconventional as the attractor moduli points can be far from the so called large complex limits where most work of mirror symmetry has been done.\\

How is this problem related to stability conditions? It turns out that due to the very special property of attractor backgrounds the central charges of special Lagrangians in $K3$ and $K3\times T^{2}$ have certain highly non-generic property. Roughly speaking if we have a right framework of defining stability conditions such that special Lagrangians are stable objects then this non-generic property would mean that there is an unexpected huge number of stability walls containing the attractor background. From the perspective of wall crossings this is certainly very unusual and intricate. \\

Unfortunately the right framework of defining stability conditions for derived Fukaya categories does not exist yet. Nevertheless we can apply  the mirror symmetry of $K3$ surfaces and ask if the above property is transformed into some statements which can be proved independently and can be viewed as the corresponding mirror statement under a correspondence of type (5).  In other words we want to check if there is a highly non-generic configuration of stability walls   for stability conditions on the derived category of coherent sheaves on the mirror $K3$ surface and if this behavior has nice correspondences (induced by mirror symmetry) with the similar behavior mentioned before for the original $K3$ surface in an attractor background.  If that is the case then we not only have some nontrivial results about mirror symmetry in attractor backgrounds but also can connect them to stability walls in the space  of stability conditions. Then the three topics in this introduction are finally on the same page. \\

In this paper we study this problem and the author believes that we have  verified that the idea in the previous paragraph works and the expectations are true. \\

In section 2 we review the attractor backgrounds of $K3$ surfaces given by Moore \citep{M}  and also some basic properties of them obtained by Shioda and Inose \citep{SI}. In section 3 we review the work of Bridgeland  \citep{Br} on stability conditions on derived categories of coherent sheaves on $K3$ surfaces. Then following the exposition of Dolgachev \citep{Do} and Huybrechts \citep{Hu1} we explain mirror symmetry of $K3$ surfaces. We discuss the non-generic property in an attractor background mentioned before in section 5. It is observed by physicists Aspinwall, Maloney and Simons \citep{AMS} and has inspired the author to write this paper. Our main theorems are formulated and proved in section 6. There are also some speculations at the end of section 6.

\section{Singular K3 Surfaces as Attractors}

Let us formulate the mathematical condition proposed by physicists that characterizes attractor backgrounds. For a Calabi-Yau threefold $X$ we fix a cohomology class (called a $charge$) \begin{equation}\gamma\in H^{3}(X,\mathrm{Z})\end{equation} Let $\tilde{\mathcal{M}}$ be the universal cover of the moduli space $\mathcal{M}$ of complex structures of $X$. Consider the family of the marked Calabi-Yau threefolds $\pi:\tilde{\mathcal{X}}\rightarrow \tilde{\mathcal{M}}$ and the pull-back $\tilde{\mathcal{L}}$ of the Hodge bundle. $\tilde{\mathcal{L}}:=R\pi_{*}\Omega^{3,0}(\tilde{\mathcal{X}})\rightarrow \tilde{\mathcal{M}}$. For $\gamma$ and a section $\Omega$ of $\tilde{\mathcal{L}}$ we define  a function called the $central\ charge$\begin{equation}Z(\gamma,\Omega):=\int_{X}\Omega\wedge\gamma\end{equation}
To simplify notations later we will just write $Z(\gamma)$ for $Z(\gamma,\Omega)$ when it is clear what the holomorphic 3-form $\Omega$ is equal to.

 Attractor backgrounds are determined by charges. We fix a charge $\gamma$. Then attractor backgrounds are nonvanishing local minima of $|Z(\gamma,\Omega)|^{2}$. It is shown in \citep{M} that this condition can be reformulated as the following $attractor\ equation$ \begin{equation}\gamma=\gamma^{3,0}+\gamma^{0,3}\in H^{3,0}(X)\oplus H^{0,3}(X)\end{equation}In other words, the attractor complex moduli must be such that the Hodge decomposition of the topological object $\gamma$ has only $(3,0)$ and $(0,3)$ parts. So the essence of the attractor mechanism is that special complex moduli (attractor backgrounds) are determined by cohomology classes (charges). From now on we can forget physical motivations and the mathematically precise formulation of the attractor equation will be our starting point.\\

Let us consider Calabi-Yau threefolds $X:=S\times T^{2}$ with varying complex moduli. Here $S$ is a $K3$ surface. We choose a basis of $H^{1}(T^{2},\mathrm{Z})$ which  can  be written as $\{dx, dy\}$ with  coordinates $x, y$ of $T^{2}$. We can choose appropriate $x, y$ such that we have the normalization $$\int_{T^{2}}dx\wedge dy=1$$ Since $H^{3}(X,\mathrm{Z})=H^{2}(S,\mathrm{Z})\otimes H^{1}(T^{2},\mathrm{Z})$, a charge can be written as \begin{equation}\gamma=p dx+q dy\end{equation} where $p,q \in H^{2}(S,\mathrm{Z})$. Let $\tau$ be the complex modulus of $T^{2}$ and $\Omega$ be a holomorphic 2-form on $S$, then \begin{equation}\Omega^{3,0}:=(dx+\tau dy)\wedge\Omega\end{equation} is a holomorphic 3-form on $X$. In \citep{M} the attractor equation has been solved and the result is \begin{equation}\tau={p\cdot q+i\sqrt{D_{p,q}}\over p^{2}}\end{equation}\begin{equation}\Omega=q-\bar{\tau}p\end{equation} where $p\cdot q$ is the intersection paring on $H^{2}(S, \mathrm{Z})$, $p^{2}:=p\cdot p$ and $D_{p,q}=p^{2}q^{2}-(p\cdot q)^{2}$. We always assume that $D_{p,q}>0$.

$H^{2}(S, \mathrm{Z})$ endowed with the intersection paring is a lattice isomorphic to \begin{equation}\Gamma:=2(-E_{8})\oplus 3U\end{equation} where $U$ is called the $hyperbolic\ plane$ and is defined to be the free group of rank two generated by isotropic vectors $\{e_{1}, e_{2}\}$ with $ e_{1}\cdot e_{2}=1$. The signature of $\Gamma$ is $(3,19)$. The choice of a lattice isomorphism $\phi:H^{2}(S, \mathrm{Z})\simeq L$ is called a $marking$. Recall the definition of Neron-Severi lattice $NS(S)$ and transcendental lattice $T_{S}$ (in the following $j: H^{2}(S,\mathrm{Z})\rightarrow  H^{2}(S,\mathrm{C})$ is the inclusion)\begin{equation}NS(S):=H^{1,1}(S,\mathrm{R})\cap Im\ j(H^{2}(S,\mathrm{Z}))\end{equation}\begin{equation}T_{S}:=NS(S)^{\perp}\cap Im\ j(H^{2}(S,\mathrm{Z})) \end{equation}For an algebraic $K3 $ surface $NS(S)=Pic(S)$. The rank $\rho(S)$ of $NS(S)$ is called the Picard number. The intersection paring on $NS(S)$ has signature $(1,\rho(S)-1)$.  As noticed in  \citep{M}, the equation (12) means that $NS(S)$ is the orthogonal complement of the lattice generated by $p,q$. So  $\rho(S)$ is 20. Note that 20 is the largest possible Picard number of a $K3$ surface.\\

Shioda and Inose have studied $K3$ surfaces with $\rho(S)=20$ (\citep{SI}). They are called $singular \ K3\ surfaces$\footnote {The word $singular$ here does NOT mean $nonsmooth$. Perhaps it is more appropriate to call them $attractive\ K3\ surfaces$ because of theorem 2.2.}\begin{thm} (\citep{SI})There is a natural one-to-one correspondence from the set of singular $K3$ surfaces to the set of equivalence classes of positive-definite even integral binary quadratic forms with respect to the  action of $SL_{2}(\mathrm{Z})$.\end{thm}

For a singular $K3$ surface the rank 2 transcendental lattice $T_{S}$ has a natural orientation by requiring that the ratio of two corresponding periods has positive imaginary part. Let $\{p,q\}$ be an oriented basis of $T_{S}$ then we can define a quadric form $Q_{p,q}$ by \begin{equation} Q_{p,q}:=\left(\begin{array}{cc}p^{2} & p\cdot q\\ p\cdot q & q^{2} \end{array}\right)\end{equation} This is the quadratic form mentioned in the theorem. Here positive-definite means $D_{p,q}>0$. The equivalence relation with respect to the action of $SL_{2}(\mathrm{Z})$  is that $Q_{1}\sim Q_{2}$ iff $Q_{1}=r^{T}Q_{2}r$ for some $r\in SL_{2}(\mathrm{Z})$. Note that $p^{2}$ and $q^{2}$ need not to be squares of integers.

Conversely for any positive-definite even integral binary quadratic form Shioda and Inose can construct a singular $K3$ surface such that the intersection matrix of the rank 2 transcendental lattice is the quadratic form.\\

So we have explained the following basic fact discovered in \citep{M}.
\begin{thm}For  $X:=S\times T^{2}$  and  a charge specified by the choice of $p,q$ above, the attractor background is determined  by (11) (12).  The  $K3$ surface $S$ is the singular $K3$ surface corresponding to $Q_{p,q}$.  \end{thm}

From now on a $K3$ surface in an attractor background is always understood as the singular $K3$ surface corresponding to some given $Q_{p,q}$.

Shioda and Inose have also proved that singular $K3$ surfaces always carry elliptic fibration structures.

\begin{thm}Every singular $K3$ surface has an elliptic fibration structure with an infinite group of sections.\end{thm}
This theorem implies that every singular $K3$ surface has an infinite group of automorphisms (this is actually a theorem in \citep{SI} and theorem 2.3 is obtained in the proof of it). It is also interesting to know that there might be more than one elliptic fibration structures with sections on a given singular $K3$ surface.

\section{Stability Conditions on K3 Surfaces}

Bridgeland has defined stability conditions on a triangulated category based on the work of Douglas. We summarize some results from \citep{Br}.

\begin{dfn} A $stability\ condition$ $\sigma=(Z,\mathcal{P})$ on a triangulated category $\mathcal{D}$ consists of a linear map $Z:K(\mathcal{D})\rightarrow \mathrm{C}$ called the $central\ charge$ (here $K(\mathcal{D})$ is the Grothendieck group of $\mathcal{D}$), and full additive subcategories $\mathcal{P}(\phi)\subset\mathcal{D}$ for each $\phi\in\mathrm{R}$, satisfying the following axioms
\begin{enumerate}\item if $0\neq E\in\mathcal{P}(\phi)$ then $Z(E)=m(E)\exp(i\pi\phi)$ for some $m(E)\in\mathrm{R}_{>0}$,\item for all $\phi\in\mathrm{R}$, $\mathcal{P}(\phi+1)=\mathcal{P}(\phi)[1]$,\item if $\phi_{1}>\phi_{2}$ and $A_{j}\in\mathcal{P}(\phi_{j})$ then $Hom_{\mathcal{D}}(A_{1},A_{2})=0$,\item for $0\neq E\in \mathcal{D}$ there is a finite sequence of real numbers $$\phi_{1}>\phi_{2}>\cdots>\phi_{n}$$ and a collection of triangles \begin{diagram}
0=E_{0} & &\rTo & & E_{1} &  &\rTo &       & E_{2}  &  \rTo & \cdots &\rTo & E_{n-1} &          &\rTo &      &E_{n}=E\\
&\luDashto &    & \ldTo & &\luDashto &    & \ldTo & &        &        &    &         &\luDashto &    & \ldTo &\\
&          & A_{1} &    &  &         & A_{2}&     &  &       &        &    &         &          & A_{n} &    &
\end{diagram}with $A_{j}\in\mathcal{P}(\phi_{j})$ for all $j$.\end{enumerate}\end{dfn}

Each subcategory $\mathcal{P}(\phi)$  is abelian and nonzero objects of it are called $semistable\ objects$ with phase $\phi$ in $\sigma$. Simple objects of $\mathcal{P}(\phi)$ are called $stable\ objects$. Define \begin{equation}\phi_{\sigma}^{+}(E)=\phi_{1},\ \ \phi_{\sigma}^{-}(E)=\phi_{n}\end{equation}
The set of stablity conditions\footnote{In fact one imposes a technical condition called local finiteness. For details see \citep{Br}.} on a triangulated category $\mathcal{D}$ is denoted by $Stab(\mathcal{D})$. It has a natural topology induced by the following generalized metric \begin{equation}d(\sigma_{1},\sigma_{2})=\sup_{0\neq E\in\mathcal{D}}\{|\phi_{\sigma_{2}}^{-}(E)-\phi_{\sigma_{1}}^{-}(E)|,\ |\phi_{\sigma_{2}}^{+}(E)-\phi_{\sigma_{1}}^{+}(E)|,\ |\log{m_{\sigma_{2}}(E)\over m_{\sigma_{1}}(E)}|\}\end{equation}

We will not discuss abstract triangulated categories. A triangulated category in this paper is always the bounded derived category of coherent sheaves $\mathcal{D}(X)$ on a smooth complex projective variety $X$. For two objects $E,F\in \mathcal{D}(X)$ we define the Euler form $\chi(E,F)$ by $$\chi(E,F):=\sum_{i}(-1)^{i} \dim_{\mathrm{C}}Hom_{X}^{i}(E,F[i])$$ Define $\mathcal{N}(X)=K(X)/K(X)^{\perp}$ where $K(X)^{\perp}$ is the radical of the bilinear form on $K(X)$ induced by the Euler form.

From now on we shall always assume that a stability condition $\sigma=(Z,\mathcal{P})$ is $numerical$ in the sense that the central charge $Z$ of it is  of the following form \begin{equation}Z(E)=-\chi(\pi(\sigma),E)\end{equation} for some vector $\pi(\sigma)\in\mathcal{N}(X)\otimes\mathrm{C}$. The set of locally finite numerical stability conditions on $\mathcal{D}(X)$ is denoted by $Stab(X)$ instead of $Stab(\mathcal{D})$. From now on whenever we say space of stability conditions we always mean space of locally finite numerical stability conditions.

\begin{thm} For each connected component of $Stab^{*}(X)\in Stab(X)$ there is a linear subspace $V\in\mathcal{N}(X)\otimes\mathrm{C}$ such that $$\pi: Stab^{*}(X)\rightarrow\mathcal{N}(X)\otimes\mathrm{C} $$ is a local homeomorphism onto an open subset of $V$.\end{thm}

Now let $S$ be a projective $K3$ surface. One can show that\begin{equation}\mathcal{N}(S)=\mathrm{Z}\oplus NS(S)\oplus\mathrm{Z}\end{equation}On the cohomology ring $H^{*}(S, \mathrm{Z})=H^{0}(S, \mathrm{Z})\oplus H^{2}(S, \mathrm{Z})\oplus H^{4}(S, \mathrm{Z})$ we define the $Mukai \ pairing$ $(\cdot,\cdot)$\begin{equation}((r_{1},D_{1},s_{1}),(r_{2},D_{2},s_{2})):=D_{1}\cdot D_{2}-r_{1}s_{2}-r_{2}s_{1}\end{equation} For $E\in \mathcal{D}(S)$, the $Mukai\ vector$ $v(E)$ is \begin{equation}v(E)=(r(E),c_{1}(E),s(E)):=ch(E)\sqrt{td(S)}\end{equation} One can show that $v(E)$ is integral and lies in $\mathcal{N}(S)$. Moreover the condition (19) becomes \begin{equation}Z(E)=(\pi(\sigma),v(E))\end{equation}for some vector $\pi(\sigma)\in\mathcal{N}(S)\otimes\mathrm{C}$.\\

To study the space of stability conditions on $S$ we first need to find some stability conditions. Take a pair of divisors $B,\omega\in NS(S)\otimes\mathrm{R}$ such that $\omega$ lies in the ample cone. Define a group homomorphism $Z: \mathcal{N}(S)\rightarrow \mathrm{C}$ by \begin{equation}Z(E):=(\exp(B+i\omega), v(E))\end{equation}If we write it explicitly it is \begin{equation}Z(E)={1\over 2r}((D^{2}-2rs)+r^{2}\omega^{2}-(D-rB)^{2})+i(D-rB)\cdot \omega\end{equation}where $(r,D,s):=v(E)$ and we assume $r\neq0$. If $r=0$ then\begin{equation}Z(E)=(D\cdot B-s)+i(D\cdot\omega)\end{equation}

\begin{thm}Suppose $B,\omega\in NS(S)\otimes\mathrm{Q}$ and $\omega$ lies in the ample cone. Suppose $Z(E)$ does not lie in $\mathrm{R}_{\leq0}$ for all spherical sheaves $E$ (this holds if $\omega^{2}>2$). Then there is a stability condition on $\mathcal{D}(S)$ with central charge given by (24).\end{thm}

Next we describe the space of stability conditions. Let $\mathcal{P}(S)\subset\mathcal{N}(S)\otimes\mathrm{C}$ be the open set consisting of all elements of $\mathcal{N}(S)\otimes\mathrm{C}$ whose real and imaginary parts span positive definite two-planes in $\mathcal{N}(S)\otimes\mathrm{R}$. It has two connected components exchanged by complex conjugation. Consider  the following  tube domain $$\{B+i\omega\in NS(S)\otimes\mathrm{C}:\omega^{2}>0\}$$ It is easy to see that if $B+i\omega$ is in this tube domain then\begin{equation}\Psi:=\exp(B+i\omega) \in \mathcal{P}(S)\end{equation}So $\mathcal{P}(S)$ contains the space of complexified Kahler deformations (the tube domain).

Let $\mathcal{P}^{+}(S)$ be the connected component of $\mathcal{P}(S)$ containing the vectors of the form $\exp(B+i\omega)$ for ample real classes $\omega$. Define $$\Delta(S):=\{\delta\in\mathcal{N}(S): (\delta,\delta)=-2\}$$ Note that it contains $(-2)$-classes of curves of $S$. For each $\delta\in\Delta(S)$, we define a hyperplane $$\delta^{\perp}:=\{\Psi\in\mathcal{N}(S)\otimes\mathrm{C}:(\Psi,\delta)=0\}$$ The hyperplane complement is $$\mathcal{P}_{0}(S)=\mathcal{P}(S)\setminus \bigcup_{\delta\in\Delta(S)}\delta^{\perp}$$ We also define $$\mathcal{P}^{+}_{0}(S)=\mathcal{P}^{+}(S)\setminus \bigcup_{\delta\in\Delta(S)}\delta^{\perp}$$

 A connected component $Stab^{*}(S)\subset Stab(S)$ is a $good\ component$ if it contains a point $\sigma$ such that $\pi(\sigma)\in \mathcal{P}_{0}(S)$. Let $Stab^{0}(S)\subset Stab(S)$ be the connected component containing the stability conditions constructed in theorem 3.2.
\begin{thm} The subset $\mathcal{P}_{0}(S)\subset \mathcal{N}(S)\otimes\mathrm{C}$ is open and the restriction\begin{equation}\pi:\pi^{-1}(\mathcal{P}_{0}(S))\rightarrow \mathcal{P}_{0}(S)\end{equation} is a topological covering map. If a connected component $Stab^{*}(S)$ is good then the image $\pi(Stab^{*}(S))$ contains one of the two connected components of the open subset $\mathcal{P}_{0}(S)$. Moreover\begin{equation}\pi:Stab^{0}(S)\rightarrow \mathcal{P}^{+}_{0}(S)\end{equation}  is a topological covering map. \end{thm}

Note that even though in theorem 3.2 we have assumed that $B$ and $\omega$ are rational  we can drop this restriction now as along as $\exp(B+i\omega)\in \mathcal{P}_{0}(S)$. Theorem 3.3 guarantees that there are stability conditions in the pre-image  under the map $\pi$ and  central charges are still given by (24).\\

There are certain walls in the space of stability conditions that are our main interest in this paper. We say a set of objects $T\subset\mathcal{D}(S)$ has bounded mass in a connected component $Stab^{*}(S)$ if $\sup\{m_{\sigma}(E): E\in T\}<\infty$ for some $\sigma\in Stab^{*}(S)$.

\begin{thm}Suppose the subset $T$ has bounded mass in a good component $Stab^{*}(S)$ and fix a compact subset $B\subset Stab^{*}(S)$. Then there is a finite collection $\{W_{\tilde{\gamma}}:\tilde{\gamma}\in\tilde{\Gamma}\}$ of real codimension one submanifolds of $Stab^{*}(S)$ such that any connected component $$C\subset B\setminus\bigcup_{\tilde{\gamma}\in\tilde{\Gamma}}W_{\tilde{\gamma}}$$ has the following property: if $E\in T$ is semistable in $\sigma$ for some $\sigma\in C$ the $E$ is semistable  for all $\sigma\in C$. \end{thm}

Let us call the real codimension one submanifolds $W_{\tilde{\gamma}}$ $stability\ walls$. The above theorem means that when we change the stability condition the set of semistable objects does not change if no stability walls are crossed. If a stability wall is crossed then the set of semistable objects might change and determining these changes  quantitatively is a subject called wall crossing formulas.

These stability walls can be described explicitly. In fact let $\{v_{i}:i\in I\}$ be the finite set of Mukai vectors of objects of set of nonzero objects $A\in \mathcal{D}(S)$ such that for some $\sigma\in B$ and some $E\in T$ one has $m_{\sigma}(A)\leq m_{\sigma}(E)$. Let $\mathrm{\tilde{\Gamma}}$ be the set of pairs $i,j\in I$ such that $v_{i}$ and $v_{j}$ do not lie on the same real line in $\mathcal{N}(S)\otimes\mathrm{R}$. Then the definition in \citep{Br} is
\begin{dfn} (Definition-Theorem) The stability walls in theorem 3.4  are defined by \begin{equation}W_{\tilde{\gamma}}=\{\sigma=(Z,\mathcal{P})\in Stab^{*}(S):Z(v_{i})/Z(v_{j})\in \mathrm{R}_{>0}\}\end{equation}\end{dfn}

\section{Mirror Symmetry of K3 Surfaces}

We will be following the expositions of \citep{Hu1} and \citep{Do}.\\

First we introduce some period domains. Let $\Gamma$ be the $K3$ lattice $2(-E_{8})\oplus3U$ and $\Gamma_{\mathrm{R}}=\Gamma\otimes\mathrm{R}$. Define $Gr_{2}^{po}(\Gamma_{\mathrm{R}})$ to be the space of all oriented two dimensional positive subspaces of $\Gamma_{\mathrm{R}}$. Here positive means that the restriction of the bilinear form on $\Gamma_{\mathrm{R}}$ induced by the lattice paring is positive definite. Similarly we can define $Gr_{3}^{po}(\Gamma_{\mathrm{R}})$ and $Gr_{4}^{po}(\Gamma_{\mathrm{R}}\oplus U_{\mathrm{R}})$.  We also define\begin{equation}Gr_{2,1}^{po}(\Gamma_{\mathrm{R}}):=\{(P,\omega): P\in Gr_{2}^{po}(\Gamma_{\mathrm{R}}),\omega\in P^{\perp}\subset \Gamma_{\mathrm{R}},\omega^{2}>0\}\end{equation}and \begin{equation}Gr_{2,2}^{po}(\Gamma_{\mathrm{R}}\oplus U_{\mathrm{R}}):=\{(H_{1},H_{2}): H_{i}\in Gr_{2}^{po}(\Gamma_{\mathrm{R}}\oplus U_{\mathrm{R}}),H_{1}\perp H_{2}\}\end{equation}

Take a copy of hyperbolic plane $U$ (which is not to be considered as a sublattice of $\Gamma$). Let us fix a standard basis $\{w,w^{*}\}$ of $U$ (so $w^{2}=(w^{'})^{2}=0$ and $w\cdot w^{'}=1$). There is a natural injection $s: Gr_{2,1}^{po}(\Gamma_{\mathrm{R}})\times\Gamma_{\mathrm{R}}\rightarrow Gr_{2,2}^{po}(\Gamma_{\mathrm{R}}\oplus U_{\mathrm{R}})$. \begin{equation}s(((P,\omega),B))= (H_{1},H_{2})\end{equation}  with $$H_{1}:=\{x-(x\cdot B)w:x\in P\}$$ $$H_{2}:= (1/2(\omega^{2}-B^{2})w+w^{*}+B)\mathrm{R}\oplus(\omega-(\omega\cdot B)w)\mathrm{R}$$
There is an isomorphism \begin{equation}\varphi: Gr_{3}^{po}(\Gamma_{\mathrm{R}})\times \mathrm{R}_{>0}\times\Gamma_{\mathrm{R}}\cong Gr_{4}^{po}(\Gamma_{\mathrm{R}}\oplus U_{\mathrm{R}}) \end{equation}It is given by\begin{equation}\varphi: (F,\alpha,B)\rightarrow \Pi:=B^{'}\mathrm{R}\oplus F^{'}\end{equation}where $F^{'}:=\{f-(f\cdot B)w: f\in F\}$ and $B^{'}:=B+1/2(\alpha-B^{2})w+w^{*}$. The inverse map $\psi: \Pi\rightarrow (F,(B^{'})^{2},B)$ is defined in the following way. First define a three dimensional space $F^{'}:=\Pi\cap w^{\perp}$ and then define $F:=\pi(F^{'})$ where $\pi$ is the natural projection $\pi:\Gamma_{\mathrm{R}}\oplus U_{\mathrm{R}}\rightarrow \Gamma_{\mathrm{R}}$. There exists a $B^{'}\in \Pi$ such that $\Pi=B^{'}\mathrm{R}\oplus F^{'}$ is an orthogonal splitting. Such a $B^{'}$ is of course not unique but we can make it unique by imposing the condition $B^{'}\cdot w=1$. Finally we define $B:= \pi(B^{'})$. We call $B$ a $B-field$.\\

The spaces defined above  have interpretations as period domains.
\begin{dfn} A $marked\ K3\ surface$ is a pair $(S,\phi)$ such that $S$ is an $K3$ surface and $\phi: H^{2}(S,\mathrm{Z})\cong\Gamma$ is an lattice isomorphism (it must respect the parings). Two marked $K3$ surfaces are equivalent if there is an isomorphism of the underlying $K3$ surfaces such the lattice isomorphisms are identical after composing with the induced map on the cohomology ring. \end{dfn}

Then we can define the period map $P^{complex}$ by sending an equivalence class of marked $K3$ surfaces $(S,\phi)$ to \begin{equation}P^{complex}(S,\phi):=\phi\langle Re \Omega, Im \Omega\rangle \in Gr_{2}^{po}(\Gamma_{\mathrm{R}})\end{equation}
Here $\Omega$ is an holomorphic 2-form on $S$ and $\langle Re \Omega, Im \Omega\rangle$ is the subspace generated by $Re \Omega$ and $Im \Omega$.

We can also define  a marked\ complex\ hyperkahler\ $K3$ surface with\ a\ $B$-field which is a tuple $(S,g,I,B,\phi)$,  a marked hyperkahler $K3$ surface with a $B$-field which is a tuple $(S,g,B,\phi)$ and a marked Kahler $K3$ surface with a $B$-field which is a tuple $(S,\omega,B,\phi)$. The meanings of $S$ and $\phi$ are the same as before. $g$ is a hyperkahler metric on $S$, $I$ is a compatible complex structure and $B\in H^{2}(S,\mathrm{R})$. $\omega$ is a Kahler class on $S$. We can also define the obvious equivalence relations. By Yau's solution of Calabi's conjecture we know that there is a natural bijection between  the set of  equivalence classes of marked\ complex\ hyperkahler\ $K3$ surfaces with\  $B$-fields and the set of equivalence classes of  marked Kahler $K3$ surfaces with $B$-fields.

We have period maps from the set of equivalence classes of $(S,g,B,\phi)$ and respectively $(S,g,I,B,\phi)$ (or equivalently $(S,\omega,B,\phi)$) to period domains.

Let $H^{2}_{+}(S, \mathrm{R})$ be the three dimensional space spanned by  $Re \Omega$, $Im \Omega$ and $\omega$. It is the cohomology of the space of self-dual 2-forms (see \citep{A1}) and is determined by the metric $g$. Now we define the period maps $$P^{4,4}(S,g,B,\phi):=$$\begin{equation}=(\phi(H^{2}_{+}(X,\mathrm{R})), \omega^{2}, \phi(B))\in Gr_{3}^{po}(\Gamma_{\mathrm{R}})\times \mathrm{R}_{>0}\times\Gamma_{\mathrm{R}}\cong Gr_{4}^{po}(\Gamma_{\mathrm{R}}\oplus U_{\mathrm{R}})\end{equation}

\begin{equation}P^{2,2}(S,g,I,B,\phi):=(P(S,\alpha,\phi),\phi(B))\in Gr_{2,1}^{po}(\Gamma_{\mathrm{R}})\times \Gamma_{\mathrm{R}}\end{equation}Here we have used the fact that there is a natural bijection between the set of equivalence classes of $(S,g,I,\phi)$ and the set of equivalence classes of $(S,\alpha,\phi)$.   $P(S,\alpha,\phi)$ is defined to be \begin{equation}P(S,\alpha,\phi):=(P^{complex}(S,\phi),\phi(\alpha))\in Gr_{2,1}^{po}(\Gamma_{\mathrm{R}})\end{equation}
These period maps are all $O(\Gamma)$-equivariant.

We recall the following famous surjectivity theorem \citep{CCS}.
\begin{thm}$P^{complex}$ defined in (36) is surjective.\end{thm}

Let us assume that there is a splitting $\Gamma=\Gamma^{'}\oplus U^{'}$ where $U^{'}$ is a copy of hyperbolic plane. We fix an isomorphism $\mathrm{f}:U^{'}\cong U$. Let $\iota$ be the involution given by $\iota: (H_{1},H_{2})\rightarrow (H_{2},H_{1})$. Let $\xi\in O(\Gamma\oplus U)$ to be the map which is the identity on $\Gamma^{'}$ and swaps $U$ and $U^{'}$ via the isomorphism $\mathrm{f}$.
\begin{dfn} Define the $mirror\ map$ on the period domain $\tilde{\xi}$  to be \begin{equation}\tilde{\xi}:=\iota\circ\xi\end{equation} It acts on $Gr_{2,2}^{po}(\Gamma_{\mathrm{R}}\oplus U_{\mathrm{R}})$ and is an involution.\end{dfn}
Since there is an natural projection $\pi:Gr_{2,2}^{po}(\Gamma_{\mathrm{R}}\oplus U_{\mathrm{R}})\rightarrow Gr_{4}^{po}(\Gamma_{\mathrm{R}}\oplus U_{\mathrm{R}})$ by sending $(H_{1},H_{2})$ to $H_{1}\oplus H_{2}$, $\tilde{\xi}$ also naturally acts on $Gr_{4}^{po}(\Gamma_{\mathrm{R}}\oplus U_{\mathrm{R}})$ and coincides with $\xi$. We also recall that $Gr_{2,1}^{po}(\Gamma_{\mathrm{R}})\times \Gamma_{\mathrm{R}}$ can be considered as a subspace of $Gr_{2,2}^{po}(\Gamma_{\mathrm{R}}\oplus U_{\mathrm{R}})$ via a natural injection $s$.

\begin{thm}Let $\{v,v^{*}\}$ be a basis of $U^{'}$ that corresponds to $\{w,w^{*}\}$ of $U$ under $\mathrm{f}$. Let $pr$  be the natural projection $pr:\Gamma_{\mathrm{R}}\rightarrow\Gamma^{'}_{\mathrm{R}}$. Let $((P,\omega),B)\in Gr_{2,1}^{po}(\Gamma_{\mathrm{R}})\times \Gamma_{\mathrm{R}}$ such that $\omega, B\in\Gamma^{'}_{\mathrm{R}}\oplus\mathrm{R}v$. Then its image under the mirror map $\tilde{\xi}$ is $((\check{P},\check{\omega}),\check{B})\in Gr_{2,1}^{po}(\Gamma_{\mathrm{R}})\times \Gamma_{\mathrm{R}}$
\begin{equation}\check{\Omega}:={1\over (Re \Omega)\cdot v}(pr(B+i\omega)-1/2(B+i\omega)^{2}v+v^{*})\end{equation}\begin{equation}\check{B}+i\check{\omega}:={1\over (Re\Omega)\cdot v}(pr(\Omega)-(\Omega\cdot B)v)\end{equation} Here $P$ is spanned by $Re \Omega$ and $Im \Omega$ and we choose $\Omega$ such that $Im \Omega$ is orthogonal to $v$. \end{thm}

So far the mirror map is studied only on period domains. To relate it to actual $K3$ surfaces we use the construction of Dolgachev \citep{Do}.

\begin{dfn} Let $N$ be a sublattice of $\Gamma$ of signature $(1,r)$ ($r$ can be 0). A marked $N$-polarized  $K3$ surface is a marked $K3$ surface $(S,\phi)$ such that $N\subset\phi(Pic(S))$.\end{dfn}

For a projective $K3$ surface $S$ we define \begin{equation}\tilde{\Delta}(S):=\{\delta\in Pic(S):\delta^{2}=-2\}\end{equation} $\tilde{\Delta}(S)$ has two components \begin{equation}\tilde{\Delta}(S)=\tilde{\Delta}(S)^{+}\coprod\tilde{\Delta}(S)^{-}\end{equation}with $\tilde{\Delta}(S)^{+}= -\tilde{\Delta}(S)^{-}$ and $\tilde{\Delta}(S)^{+}$ consists of effective classes. Define\begin{equation}V(S):=\{x\in H^{1,1}(S)\cap H^{2}(S,\mathrm{R}): x^{2}>0\}\end{equation} This is a cone consisting of two components.  We denote by $V(S)^{+}$ the component which contains the class of some Kahler form on $S$.
Then we define \begin{equation}C(S):=\{x\in V(S)^{+}:  (x,\delta)\geq 0 \ for\ any\ \delta\in \tilde{\Delta}(S)^{+}\}\end{equation}One can show that the set of interior points $C(S)^{+}$ of $C(S)$ is the Kahler cone of $S$. We also define \begin{equation}Pic(S)^{+}:=C(S)\cap H^{2}(S,\mathrm{Z}),\ \ Pic(S)^{++}:=C(S)^{+}\cap H^{2}(S,\mathrm{Z}) \end{equation}
For the sublattice $N$ we define \begin{equation}V(N):=\{x\in N_{\mathrm{R}}: x^{2}>0\}\end{equation} It has two components. We fix one  of them and call it $V(N)^{+}$. Let\begin{equation}\tilde{\Delta}(N):=\{\delta\in N: \delta^{2}=-2\}\end{equation} $\tilde{\Delta}(N)=\tilde{\Delta}(N)^{+}\coprod \tilde{\Delta}(N)^{-}$ with $\tilde{\Delta}(N)^{+}= -\tilde{\Delta}(N)^{-}$ and $\tilde{\Delta}(N)^{+}$ has the following property.  If $\delta_{1},\cdots,\delta_{k}\in \tilde{\Delta}(N)^{+}$ and $\delta$ is a linear combination of them with nonnegative coefficients then $\delta\in\tilde{\Delta}(N)^{+}$. Denote by $C(N)^{+}$ the following set\begin{equation}C(N)^{+}:\{h\in V(N)^{+}\cap N: h\cdot \delta>0 \ for\ any\ \delta\in\tilde{\Delta}(N)^{+}\}\end{equation}

\begin{dfn} A marked $N$-polarized $K3$ surface $S$ is called a marked pseudo-ample $N$-polarized $K3$ surface if \begin{equation}C(N)^{+}\cap \phi(Pic(S)^{+})\neq\emptyset\end{equation} It is called a marked ample $N$-polarized $K3$ surface if \begin{equation}C(N)^{+}\cap \phi(Pic(S)^{++})\neq\emptyset\end{equation}\end{dfn}

Let $\mathcal{T}_{\Gamma}^{2,2}$ be the set of equivalence classes of marked Kahler $K3$ surfaces with $B$-fields. Let $\mathcal{T}_{N\subset\Gamma}^{2,2}\subset \mathcal{T}_{\Gamma}^{2,2}$ be the subset consisting of all marked Kahler $K3$ surfaces with $B$-fields $(S,\omega,B,\phi)$ such that $\omega,B\in N_{\mathrm{R}}$ and $N\subset\phi(Pic(S))$. $\mathcal{T}_{N\subset\Gamma}^{2,2}$ is realized  as a subset of $Gr_{2,1}^{po}(\Gamma_{\mathrm{R}})\times \Gamma_{\mathrm{R}}$ by the period map. It is shown that the closure $\overline{\mathcal{T}}_{N\subset\Gamma}^{2,2}$ is  $$\overline{\mathcal{T}}_{N\subset\Gamma}^{2,2}=\{((P,\omega),B)\in Gr_{2,1}^{po}(\Gamma_{\mathrm{R}})\times \Gamma_{\mathrm{R}}:B,\omega\in N_{\mathrm{R}}, P\subset N^{\perp}_{\mathrm{R}}\}$$

We assume additionally that the orthogonal complement $N^{\perp}\subset\Gamma$ contains a hyperbolic plane $U^{'}$ (this implies that $r<19$). Then one can show that there is a splitting $N^{\perp}=\check{N}\oplus U^{'}$ for a sublattice $\check{N}$ with signature $(1,18-r)$.
Moreover we have \begin{equation}\Gamma=(U^{'})^{\perp}\oplus U^{'}\end{equation} It is easy to see that $\check{\check{N}}=N$.

Then we can introduce $\overline{\mathcal{T}}_{\check{N}\subset\Gamma}^{2,2}$ and $\mathcal{T}_{\check{N}\subset\Gamma}^{2,2}$ by replacing $N$ by $\check{N}$ in the definition of $\overline{\mathcal{T}}_{N\subset\Gamma}^{2,2}$ and $\mathcal{T}_{N\subset\Gamma}^{2,2}$.

We define the period domain of marked $N$-polarized $K3$ surfaces \begin{equation}D_{N}:=\{\mathrm{C}\Omega\in \mathrm{P}(N^{\perp}\otimes\mathrm{C}):\Omega\cdot\Omega=0,\Omega\cdot\bar{\Omega}>0\}\end{equation} Here $\mathrm{P}(N^{\perp}\otimes\mathrm{C})$ is the projective space associated to $N^{\perp}\otimes\mathrm{C}$. $D_{N}$ is a subset of \begin{equation}Q_{\Gamma}:=\{\Omega\in \mathrm{P}(\Gamma_{\mathrm{C}}):\Omega^{2}=0, \Omega\cdot\bar{\Omega}>0\} \cong Gr^{po}_{2}(\Gamma_{R})\end{equation} We define the following tube domain\begin{equation}T_{N}:=\{B+i\omega\in N\otimes\mathrm{C}:\omega^{2}>0\}\end{equation} $D_{N}$ and $T_{N}$ are interpreted as the space of complex deformations and the space of complexified Kahler deformations respectively  for a marked $N$-polarized $K3$ surface (but note that in the definition above we do not insist that $\omega$ is Kahler). Similarly we can define $D_{\check{N}}$ and $T_{\check{N}}$.

Clearly the map $((P,\omega),B)\rightarrow B+i\omega$ defines a surjection $\overline{\mathcal{T}}_{N\subset\Gamma}^{2,2}\rightarrow T_{N}$. There is also the natural surjective map $\overline{\mathcal{T}}_{N\subset\Gamma}^{2,2}\rightarrow D_{N}$ by forgetting $\omega$ and $B$.

\begin{thm}The mirror map $\tilde{\xi}$ induces a bijection \begin{equation}\overline{\mathcal{T}}_{N\subset\Gamma}^{2,2}\cong\overline{\mathcal{T}}_{\check{N}\subset\Gamma}^{2,2}\end{equation} We also have the following commutative diagram \begin{equation}\begin{diagram}
\overline{\mathcal{T}}_{N\subset\Gamma}^{2,2}& \rTo &\overline{\mathcal{T}}_{\check{N}\subset\Gamma}^{2,2}\\
\dTo & & \dTo\\
T_{N} & \rTo & D_{\check{N}}\\
\end{diagram}\end{equation}where the two vertical maps are the two surjective maps defined above. The top horizontal one is the bijection  $\overline{\mathcal{T}}_{N\subset\Gamma}^{2,2}\cong\overline{\mathcal{T}}_{\check{N}\subset\Gamma}^{2,2}$ and the bottom horizontal one is an isomorphism $a$ given by \begin{equation}a(z):=[z-{1\over 2}z^{2}v+v^{*}], z\in T_{N}\end{equation}\end{thm} Clearly the above theorem is still true if we swap $N$ and $\check{N}$ because $\check{\check{N}}=N$.

 We also define \begin{equation}\Gamma^{'}:=(U^{'})^{\perp}=N\oplus \check{N}\end{equation} so that $\Gamma=\Gamma^{'}\oplus U^{'}$.

\begin{note}It is clear that what matters is  the splitting $\Gamma=\Gamma^{'}\oplus U^{'}$, the splitting $\Gamma^{'}_{\mathrm{R}}=N_{\mathrm{R}}\oplus\check{N}_{\mathrm{R}}$ and the isomorphism $\mathrm{f}:U^{'}\cong U$. The lattices $N$ and $\check{N}$ are not important.\end{note}

We do not expect that we have a bijection $\mathcal{T}_{N\subset\Gamma}^{2,2}\cong\mathcal{T}_{N\subset\Gamma}^{2,2}$. To relate period domains to $K3$ surfaces we use the following results. As pointed out in \citep{Do}, the fine moduli space of  marked $N$-polarized $K3$ surface $\mathcal{K}_{N}$ exists and the period map  $\mathrm{p}$ (which is defined to be the restriction of the period map of  marked Kahler $K3$ surfaces) maps it to $D_{N}$.
\begin{thm}Let $p$ be the restriction of the period map $\mathrm{p}:\mathcal{K}_{N}\rightarrow D_{N}$ to the subset $\mathcal{K}^{pa}_{N}$ of equivalence  classes of marked pseudo-ample  $N$-polarized $K3$ surfaces. $p$ is surjective. There is  a natural bijection between the fiber of the map $p$ over a point in $D_{N}$ and a subgroup of isometries of $\Gamma$ generated by some reflections (for details see \citep{Do}). \end{thm}

 This theorem tells us that mirror symmetry (the mirror map) exchanges  complex deformations and complexfied Kahler deformations\footnote{Note that  we do not know if $\check{\omega}$ is a Kahler class even if we call $T_{\check{N}}$ the space of complexified Kahler deformations. } associated with marked pseudo-ample $N$-polarized $K3$ surfaces and marked pseudo-ample $\check{N}$-polarized $K3$ surfaces and vice versa.  But for a given  marked pseudo-ample $N$-polarized $K3$ surface (whose complex structure and complexified Kahler structure are given by $P$ and $B+i\omega$ in the triple $((P,\omega),B))$) there is not a unique  mirror $K3$ surface even though $((\check{P},\check{\omega}),\check{B})$ is uniquely determined. That is  because the map $p$ of theorem 4.4 is not bijective even though the fiber is discrete. The following theorem improves the situation.

For any $\delta\in\tilde{\Delta}(N^{\perp})$ define \begin{equation}H_{\delta}:=\{z\in N^{\perp}_{\mathrm{C}}:z\cdot\delta=0\},\ D_{N}^{\circ}:=D_{N}\setminus(\bigcup_{\delta\in\tilde{\Delta}(N^{\perp})}H_{\delta}\cap D_{N})\end{equation}

\begin{thm}Let $\mathcal{K}_{N}^{a}$ be the subset of $\mathcal{K}_{N}$ consisting of equivalence classes of marked ample $N$-polarized $K3$ surfaces. The restriction of the period map $p$ to $\mathcal{K}_{N}^{a}$ (also denoted by $p$) induces a bijection \begin{equation}p:\mathcal{K}_{N}^{a}\rightarrow D_{N}^{\circ} \end{equation}\end{thm}

\section{Special Langangians in Attractor Backgrounds}

Let $Y=S\times T^{2}$ where $S$ is an singular $K3$ surface. We pick a hyperkahler metric $g$ such that the underlying complex structure of the singular $K3$ surface $S$ is compatible with $g$. Such a hyperkahler metric exists abundantly  and  is determined by the choice of a Kahler class on $S$. The underlying complex structure of the singular $K3$ surface is denoted by $J$.  The set of all compatible complex structures is identified with $P^{1}$ and generated by three complex structures $I,J,K$ (i.e. it is $\{aI+bJ+cK:a^{2}+b^{2}+c^{2}=1\}$). Let $\omega_{I},\omega_{J},\omega_{K}$ be the three corresponding Kahler classes and $\Omega_{I},\Omega_{J},\Omega_{K}$ be the three corresponding holomorphic 2-forms (appropriately normalized).  Then we have \begin{equation}\Omega_{I}=\omega_{J}+i\omega_{K},\ \Omega_{J}=\omega_{K}+i\omega_{I}\end{equation}

The trick of $hyperkahler\ rotation$ is that if we consider the compatible complex structure $I$  then a holomorphic cycle $L$ in $J$ becomes a special Lagrangian cycle. In fact $\omega_{I}=Im \Omega_{J}$ vanishes on $L$ since $L$ is holomorphic in $J$. And $Im \Omega_{I}=\omega_{K}=1/2(\Omega_{J}+\bar{\Omega}_{J})$ restricts to $0$ on $L$ while $Re \Omega_{I}=\omega_{J}$ restricts to the volume form on $L$. This recovers the definition of a special Lagrangian cycle. So in the complex structure $I$ each class $l\in Pic(S)$  is the Poincare dual of a special Lagrangian cycle (but it may not be a smooth special Lagrangian submanifold).

The following theorem is implicitly contained in \citep{AMS}.

\begin{thm}  Let $l_{1},l_{2}\in Pic(S)$. Pick a compatible hyperkahler metric $g$ on $S$ (the underlying complex structure is $J$) and a flat metric $g^{'}$ on $T^{2}$. Endow $Y$ with the product metric $\hat{g}$.  Take the hyperkahler rotation on $S$ by rotating $J$ to $I$.  Endow the manifold $Y$ with the product complex structure which is the product of $I$ and the one on $T^{2}$ corresponding to $\tau$. Set $u_{i}:=l_{i}\wedge dy$. Then  $u_{i}$ is dual to a special Lagrangian cycle in $Y$ endowed with the product metric and the product complex structure specified above.  Moreover \begin{equation}\arg Z(u_{1})=\pm\arg Z(u_{2})\end{equation}In fact we have a stronger result which implies (64): \begin{equation}Z(u_{i})=Z_{K3}(l_{i}):=\omega_{J}\cdot l_{i}\in \mathrm{R}\end{equation} \end{thm}

\noindent {\bf Proof} The holomorphic volume form (up to a normalization) on $S$ is $\Omega_{J}=q-\bar{\tau}p$ (equation (12))  while the holomorphic volume form (up to a normalization) on $T^{2}$ is $dz=dx+\tau dy$. After the hyperkahler rotation on  $S$ the holomorphic volume form is $\Omega_{I}= \omega_{J}+i(q-Re(\tau)p)$. So the normalized holomorphic 3-form on the product $Y$ with the product complex structure is \begin{equation}\Omega_{Y}:=\Omega_{I}\wedge dz=(\omega_{J}+i(q-Re(\tau)p))\wedge (dx+\tau dy)\end{equation}$l_{i}$  is dual to a holomorphic cycle $L_{i}$ of $S$ in $J$. Then $u_{i}$ is dual to a three cycle $L_{i}\times S^{1}$ in $Y$.  Here $S^{1}$ is a circle in $T^{2}$ dual to $dy$.

On $Y$ with the product metric and the product complex structure we should consider the restriction of $\Omega_{Y}$. Clearly  $Im\Omega_{Y}$ restricts on $L_{i}\times S^{1}$ to $0$ while $Re\Omega_{Y}$ restricts to $\omega_{J}\wedge dx|_{L_{i}\times S^{1}}$. The integration of $\omega_{J}\wedge dx|_{L_{i}\times S^{1}}$ on $L_{i}\times S^{1}$ is its normalized volume in $\hat{g}$. So we know $L_{i}\times S^{1}$ is a special Lagrangian cycle in the product $Y$.

Direct calculation gives \begin{equation}Z(u_{i})=Z(u_{i},\Omega_{Y})=\omega_{J}\cdot l_{i}\in \mathrm{R}\end{equation}Here we have used $\int_{T^{2}}dx\wedge dy=1$ and $p\cdot l_{i}=q\cdot l_{i}=0$. $\diamondsuit$\\
\begin{note}We have defined $Z_{K3}(l_{i})$ to indicate that the central charge is actually associated to the $K3$ surface. It will be related by the mirror symmetry to another central charge $Z(\mu(l_{i}))$ associated to a $K3$ surface defined in section 6.\end{note}

In fact even if $S$ is just a $K3$ surface the hyperkahler rotation trick gives us special Lagrangian cycles assuming  $S$ has holomorphic cycles. One (trivial) interesting point about a singular $K3$ is that there are many such cycles. Theorem 5.1 gives us another interesting fact.  It means that there are many special Lagrangians with aligned or anti-aligned central charges in an attractor background of $S\times T^{2}$ (after a hyperkahler rotation on $S$).  If we fix the topological data (the cohomology classes) then we can view the central charge as a function of the complex moduli. The condition that two charges having the same phase for the central charges is a condition which gives real codimension one locus (a wall) in the moduli space. So we have just shown that many codimension one walls intersect at an attractor background which is a quite special situation.

As explained in the introduction, we expect that there exist stability conditions on derived Fukaya categories such that special Lagrangians are stable objects (see \citep{A2}\citep{T} \citep{TY}  for some hints). Then the condition of having aligned central charges will become the definition of stability walls in the space of such stability conditions. So we expect that theorem 5.1 can be interpreted as the statement that  many stability walls intersect at an attractor background.

\section{Mirror Symmetry and Stability Conditions of K3 Surfaces in Attractor Backgrounds}

Let us assume that $S$ is an elliptic $K3$ surface $\pi: S\rightarrow P^{1}$ with a section $\sigma_{0}$.  We denote the cohomology class of a general fiber by $f$ and denote the section class also by $\sigma_{0}$. Then there is a hyperbolic plane sublattice $U^{'}\subset H^{2}(S,\mathrm{Z})$ generated by the basis $\{v,v^{*}\}$ with \begin{equation}v:=f,\ \ v^{*}:=f+\sigma_{0}\end{equation}Here we have used $f^{2}=0$, $f\cdot\sigma_{0}=1$ and $\sigma^{2}_{0}=-2$.

Denote the complex structure of the elliptic $K3$ surface $S$ by $J$. Pick a compatible hyperkahler metric $g$. The set of all compatible complex structures is generated by three complex structures $I,J,K$. The trick of hyperkahler rotation then tells us  that if we consider the compatible complex structure $I$ then  $\pi: S\rightarrow P^{1}$ becomes a (singular) special Lagrangian torus fibration.

Inspired by the SYZ mirror conjecture we consider  mirror symmetry of $S$ in the complex structure $I$. Note that we already have $\omega_{I}\in \Gamma^{'}_{\mathrm{R}}$. We will discuss below how to pick an $N$-polarization for a singular $K3$ surface. Here let us assume that we already have a marked pseudo-ample $N$-polarization structure.  We choose $B\in N_{\mathrm{R}}$.  $\Omega_{I}$ is normalized by  \begin{equation}\omega_{I}^{2}=(Re \Omega_{I})^{2}=(Im \Omega_{I})^{2}\end{equation} Note that $Im \Omega_{I}\in \Gamma^{'}_{\mathrm{R}}$.

Under the assumption above the $\tilde{\xi}$ mirror $(\check{\Omega}_{I},\check{\omega}_{I},\check{B})$ of $(\Omega_{I}, \omega_{I}, B)$ is given by (according to theorem 4.2)
\begin{equation}\check{\Omega}_{I}= {1\over (Re\Omega_{I})\cdot f}(B+i\omega_{I}-{1\over 2}(B+i\omega_{I})^{2}f+f+\sigma_{0})\end{equation}
\begin{equation}\check{\omega}_{I}={1\over (Re\Omega_{I})\cdot f}(Im \Omega_{I}- ((Im \Omega_{I})\cdot B)f)
\end{equation}

\begin{equation}\check{B}= {1\over (Re\Omega_{I})\cdot f} (pr(Re\Omega_{I})-((Re\Omega_{I})\cdot B)f)\end{equation}

We now study  mirror symmetry of $K3$ surfaces in an attractor background. So we assume $S$ is a singular $K3$ surface corresponding to $Q_{p,q}$.  As before we choose a hyperkahler metric $g$ such that this complex structure is $J$ in a triple $(I,J,K)$. $S$ carries elliptic fibration structures with sections. Pick and then fix such a structure. Let $f$ be the fiber class of a general fiber and $\sigma_{0}$ be the cohomology class of a section.

To specify the mirror map $\tilde{\xi}$ on the period domains we need a hyperbolic plane $U$ (with basis $\{w,w^{*}\}$) and an isomorphism $\mathrm{f}:U^{'}\cong U$. Here  $U^{'}$ is the hyperplane generated by $v:=f$ and $v^{*}:=f+\sigma_{0}$. In section 4 the lattice $U$ is somewhat abstract. One natural way to geometrically identify $U$ is to identify it with $H^{0}(S,\mathrm{Z})\oplus H^{4}(S,\mathrm{Z})$.  We extend the intersection paring on $H^{2}(S,\mathrm{Z})\cong\Gamma$ to the Mukai paring on $H^{*}(S,\mathrm{Z}):=H^{0}(S,\mathrm{Z})\oplus H^{2}(S,\mathrm{Z})\oplus H^{4}(S,\mathrm{Z})$. So \begin{equation}\Gamma_{\mathrm{R}}\oplus U_{\mathrm{R}}\cong H^{*}(S,\mathrm{R})\end{equation}

We pick a basis $\{w,w^{*}\}$ of $U$. In term of the notation of section 2, we take  \begin{equation}w=(0,0,-1),\ \ w^{*}=(1,0,0) \end{equation} Clearly $(w, w^{*})=1$. \footnote{Here $w=(0,0,-1)$ because we want to use the Mukai pairing on $H^{*}$. If we use the standard intersection pairing (we can do that) then we should take $w=(0,0,1)$ so that $w\cdot w^{*}=1$. Using Mukai pairing as the pairing on $U$ introduces a somewhat confusing point that we want to clarify now. In section 4 the pairing on the lattice $U$ is denoted by $x\cdot y$. But now it should be denoted by $(x,y)$. So the condition $w\cdot w^{*}=1$ in section 4 now becomes $(w, w^{*})=1$.}

As for the isomorphism $\mathrm{f}: U^{'}\cong U$ we pick the one that maps $\{f,\sigma_{0}\}$ to $\{w, (1,0,1)\}$. So $\mathrm{f}$ maps  $\{v,v^{*}\}$ to $\{w,w^{*}\}$. It respects the Mukai pairing. \\

Fix any marking $\phi$ of $S$. We need an $N$-polarization which is at least pseudo-ample. \begin{thm} Assume that $\sqrt{D_{p,q}}$ (which is an invariant of  equivalence classes of $Q_{p,q}$) is an integer.  Denote by $(S,I)$ the $K3$ surface $S$ with the complex structure $I$ and by $Pic(S,I)$ its Picard lattice.\footnote{We should also use  $(S,J)$ and replace $Pic(S)$ by $Pic(S, J)$. But when there is no danger of confusions we may keep using $Pic(S)$.}  Then there is a marked pseudo-ample $N$-polarization of $((S,I),\phi)$ for some  sublattice $N\subset\phi(Pic(S,I))$ with signature $(1,r)$ ($r\leq rank(Pic(S,I))-1$). Moreover  $U^{'}\subset N^{\perp}$. \end{thm}

\noindent {\bf Proof} By our assumption about $\sqrt{D_{p,q}}$ the holomorphic 2-form $\Omega_{J}=q-\bar{\tau}p$ in the complex structure $J$ is rational. After the hyperkahler rotation to $I$ we have a rational Kahler class $\omega_{I}=Im \Omega_{J}=Im(\tau)p$. Therefore  $p^{2}\omega_{I}$ (which is integral) is in $Pic(S,I)$. This guarantees  the existence of a nonempty sublattice $N$. For definiteness we can choose $N=\langle p^{2}\omega_{I}\rangle$. This $N$-polarization is an ample $N$-polarization.

 Since $f, \sigma_{0}\in Pic(S,J)$, $f\cdot \omega_{I}=\sigma_{0}\cdot\omega_{I}=0$. Therefore at least for $N=\langle p^{2}\omega_{I}\rangle$ we have $U^{'}\subset N^{\perp}$. Of course the choice of $N$ in theorem 6.1 needs not to be unique. $\diamondsuit$\\

 We  choose $B\in N_{\mathrm{R}}$ after we have chosen $N$.  If we choose $B=0$ (we will do that later) then the choice of $N$ or $N_{\mathrm{R}}$ from theorem 6.1 is irrelevant (only the existence is required). In fact since  the splitting $\Gamma=\Gamma^{'}\oplus U^{'}$ is  independent of the choice of $N$ or $N_{\mathrm{R}}$ so are equations (70) (71) (72)   when $B=0$.

  Let $(\check{S},\check{\Omega}_{I},\check{\omega}_{I},\check{B},\check{\phi})$ be a mirror of $(S,\Omega_{I}, \omega_{I}, B,\phi)$. $(\check{S},\check{\phi})$ can be assumed to be a marked pseudo-ample $\check{N}$-polarized $K3$ surface by theorem 4.4. Note that although $\check{S}$ is projective it is not guaranteed that $\check{\omega}_{I}$  Kahler. Also note that $(\check{S},\check{\phi})$ may not be uniquely determined by $(S,\Omega_{I}, \omega_{I}, B,\phi)$ if the period point of $\check{\Omega}_{I}$ is not in $D^{0}_{\check{N}}$. Of course $(\check{\Omega}_{I},\check{\omega}_{I},\check{B})$ is uniquely determined by $(S,\Omega_{I}, \omega_{I}, B,\phi)$. This guarantees that results in the rest of this section hold regardless of which $(\check{S},\check{\phi})$ we choose from theorem 4.4 .

Recall the definition of the mirror map $\tilde{\xi}=\iota\circ\xi$ where $\xi$ exchanges $U$ and $U^{'}$ via $U\cong U^{'}$ and restricts to the identity on $\Gamma_{\mathrm{R}}^{'}$. Clearly we have a lattice isomorphism between two sublattices of $H^{2}(S,\mathrm{Z})$ and $H^{2}(\check{S},\mathrm{Z})$ \begin{equation}j:=\check{\phi}^{-1}\circ\xi \circ\phi:\phi^{-1}(\Gamma^{'})\rightarrow \check{\phi}^{-1}(\Gamma^{'})\end{equation}

\begin{dfn}For a cohomology class $l\in H^{2}(S,\mathrm{Z})$ we define its $mirror\ class$ $\mu(l)\in H^{*}(\check{S},\mathrm{Z})$ in the following way. If $l\in\phi^{-1}(\Gamma^{'})$ then \begin{equation}\mu(l):=j(l)\end{equation}If $l\in \phi^{-1}(U^{'})$ then \begin{equation}\mu(l):= \check{\phi}^{-1}\circ\mathrm{f}\circ\phi(l)\end{equation} Then we can extend $\mu$ to  be defined on $H^{2}(S,\mathrm{Z})\simeq \phi^{-1}(\Gamma^{'})\oplus\phi^{-1}(U^{'})$. Here we have extended $\check{\phi}$ to the isomorphism $H^{*}(\check{S},\mathrm{Z})\rightarrow \Gamma\oplus U$ by identifying $H^{0}\oplus H^{4}$ with $U$ with the specification of the basis $\{w,w^{*}\}$.  \end{dfn}\begin{note}This transformation of cohomology classes is not only demanded by the mirror symmetry on the level of period domains (and the non-explicit construction of mirrors by the surjectivity of period maps)   but also  compatible with some other  approaches  to mirror symmetry of $K3$ surfaces (for example see \citep{Ho} \citep{Mo}). \end{note}  However we are not able to describe the mirror classes explicitly as $\check{\phi}$ is obtained by  using the surjectivity  of period maps which is not an explicit construction (at least to the author). Nevertheless since $\phi,\check{\phi}$ are  lattice isomorphisms we can compute the intersection (or Mukai) pairings of classes on $S$ or $\check{S}$ by going to $\Gamma\oplus U$ and therefore omit the markings in the formulas below.\\

Now let us formulate and prove main theorems of this paper.
\begin{thm} Let $S$ be the singular $K3$ surface corresponding to $Q_{p,q}$ with a specified elliptic fibration structure with sections. Assume that $\sqrt{D_{p,q}}$ is an integer. Fix a tuple  $(S,\Omega_{I}, \omega_{I}, B,\phi)$ described above. Fix a  marked pseudo-ample $N$-polarization for some $N$ from theorem 6.1. Fix  the lattices $U$ and $U^{'}$ with basis specified above. Fix the splittings $\Gamma=\Gamma^{'}\oplus U^{'}$. Fix the isomorphism $\mathrm{f}$ specified above.

For any $l_{1}, l_{2}\in Pic(S)$ define $Z(\mu(l_{i})):=(\exp(\check{B}+i\check{\omega}_{I}),\mu(l_{i}))$ (here $\mu(l_{i})$ is considered as an element  of $\Gamma$ via  the marking $\check{\phi}$ but we have omitted the makings in our notations). Then  \begin{equation}\arg Z(\mu(l_{1}))=\pm \arg Z(\mu(l_{2}))\end{equation}In fact we have a stronger result that implies (78)\begin{equation}Z(\mu(l_{i}))\in \mathrm{R}\end{equation}\end{thm}

\noindent{\bf Proof} The holomorphic 2-form (up to a normalization) on $S$ is given by : $\Omega_{J}=q-\bar{\tau}p$. Therefore \begin{equation}\omega_{I}=Im(\Omega_{J})=Im(\tau)p\end{equation}\begin{equation}Im(\Omega_{I})={1\over2}(\Omega_{J}+\bar{\Omega}_{J})= q- Re(\tau)p\end{equation}\begin{equation}Re(\Omega_{I})= \omega_{J}\end{equation} $(\check{\Omega}_{I},\check{\omega}_{I},\check{B})$ are given by \begin{equation}\check{\Omega}_{I}= {1\over (\omega_{J})\cdot f}(B+iIm(\tau)p-{1\over 2}(B+iIm(\tau)p)^{2}f+f+\sigma_{0})\end{equation} \begin{equation}\check{\omega}_{I}={1\over (\omega_{J})\cdot f}(q- Re(\tau)p-((q- Re(\tau)p)\cdot B)f)\end{equation} \begin{equation}\check{B}={1\over (\omega_{J})\cdot f}(pr(\omega_{J})-(\omega_{J}\cdot B)f)\end{equation}where $pr$ is the projection $pr:\Gamma_{\mathrm{R}}\rightarrow\Gamma^{'}_{\mathrm{R}}$.

Since $\Gamma=\Gamma^{'}\oplus U^{'}$ we know that $l_{i}=l_{i}^{a}+l_{i}^{b}$ with $l_{i}^{a}\in\Gamma^{'}$ and $l_{i}^{b}\in U^{'}$. Since $U^{'}\subset Pic(S)=Pic(S,J)$ we see that $l_{i}^{a}, l_{i}^{b}\in Pic(S)$. Therefore to prove the theorem it suffices to prove (79) for any $l_{i}\in \Gamma^{'}$ and for any $l_{i}\in U^{'}$.

If $l_{i}\in \Gamma^{'}$ then $$Z(\mu(l_{i}))= (\exp(\check{B}+i\check{\omega}_{I}), (0,l_{i},0))$$Here we are computing the pairing in  $\Gamma\oplus U$ as explained before. So\begin{equation}Z(\mu(l_{i}))= l_{i}\cdot \check{B}+il_{i}\cdot\check{\omega}_{I}\end{equation}Since $\Omega_{J}\cdot l_{i}=0$ ($l_{i}\in Pic(S)$) and $l_{i}$ is real we get $p\cdot l_{i}=q\cdot l_{i}=0$. Also note that $f\cdot l_{i}=0$ (because $l_{i}\in \Gamma^{'}=(U^{'})^{\perp}$). Therefore $l_{i}\cdot\check{\omega}_{I}=0$ i.e. $Z(\mu(l_{i}))\in\mathrm{R}$.

For the case of $l_{i}\in U^{'}$ it suffices to show that $Z(\mu(f)), Z(\mu(\sigma_{0}))\in\mathrm{R}$. By (25) (26),
\begin{equation}Z(\mu(f))=Z((0,0,-1))=1\end{equation}\begin{equation}Im(Z(\mu(\sigma_{0})))=Im(Z((1,0,1)))=-\check{B}\cdot\check{\omega}_{I}
\end{equation} We have $pr(\omega_{J})\cdot f=0$, $f^{2}=0$ and $f\cdot p=f\cdot q=\sigma_{0}\cdot p=\sigma_{0}\cdot q=0$. We also have $\omega_{J}\cdot \Omega_{J}=0$ which implies $$\omega_{J}\cdot p=\omega_{J}\cdot q=0$$ So we see that $\check{B}\cdot\check{\omega}_{I}=0$ which means that $Z(\mu(\sigma_{0}))$ is also real. $\diamondsuit$\\

We want to show that $\exp(\check{B}+i\check{\omega}_{I})\in \mathcal{P}_{0}(\check{S})$. In general this is hard to verify and may not be true. But we have some freedom when we choose the data $(\Omega_{I}, \omega_{I}, B)$. So the idea is that maybe we can use the freedom to make sure that $\exp(\check{B}+i\check{\omega}_{I})\in \mathcal{P}(\check{S})$ is away from $\bigcup_{\delta\in\Delta(\check{S})} \delta^{\perp}$. We certainly do not want to change the singular $K3$ surface $(S,J)$. So we fix $\Omega_{J}$.  After choosing $B$  equation (83) tells us that $\check{\Omega}_{I}$ is fixed up to a real scaling factor. Now $\Delta(\check{S})$ depends on the complex structure so we had better fix the period point $\check{\Omega}_{I}$. That means after fixing $B$ once and for all the only freedom is changing $\omega_{J}$ i.e. changing the hyperkahler metric $g$. This is natural as the attractor mechanism says nothing about $\omega_{J}$. Changing $\omega_{J}$ will change $\check{B}$.

\begin{thm} We make  the assumptions of theorem 6.2. When we choose $B$ we choose\footnote{Assumptions made in this theorem may not be optimal.}  $B=0$. In addition we assume that $D_{p,q}\neq 2p^{2}$. Then there exists some Kahler class $\omega_{J}$ on $(S,J)$ such that \begin{equation}\exp(\check{B}+i\check{\omega}_{I})\in \mathcal{P}_{0}(\check{S})\end{equation}\end{thm}

\noindent {\bf Proof} With $B=0$ equations (83) (84) (85) simplify to\begin{equation}\check{\Omega}_{I}= {1\over (\omega_{J})\cdot f}(iIm(\tau)p+{1\over 2}(Im(\tau))^{2}p^{2}f+f+\sigma_{0})\end{equation} \begin{equation}\check{\omega}_{I}={1\over (\omega_{J})\cdot f}(q- Re(\tau)p)\end{equation} \begin{equation}\check{B}={1\over (\omega_{J})\cdot f}(pr(\omega_{J}))\end{equation}

We need to check that for any $\delta\in \mathcal{N}(\check{S})$ with $(\delta,\delta)=-2$ we have $(\exp(\check{B}+i\check{\omega}_{I}),\delta)\neq 0$. Here again we are computing the pairing on $\Gamma\oplus U$ and have omitted markings in the notations. Such a $\delta=(r,D,s)$ must belong to one of the following types. \begin{enumerate}\item $D^{2}=-2, r=0,s\neq0$. In this case \begin{equation}Re (\exp(\check{B}+i\check{\omega}_{I}),\delta)=D\cdot\check{B}-s\end{equation}\item $D^{2}-2rs=-2,r\neq0,s\neq0$. In this case \begin{equation}Re (\exp(\check{B}+i\check{\omega}_{I}),\delta)={1\over 2r}(-2rs+r^{2}\check{\omega}_{I}^{2}-r^{2}\check{B}^{2}+2rD\cdot \check{B})\end{equation}\item $D^{2}=-2, r\neq0,s=0$. In this case \begin{equation}Re (\exp(\check{B}+i\check{\omega}_{I}),\delta)={1\over 2r}(r^{2}\check{\omega}_{I}^{2}-r^{2}\check{B}^{2}+2rD\cdot \check{B})\end{equation}\item $\delta=(0,D,0)$, $D^{2}=-2$. In this case \begin{equation}Re (\exp(\check{B}+i\check{\omega}_{I}),\delta)= D\cdot\check{B}\end{equation}\end{enumerate}

Pick a reference Kahler class $\omega_{J}^{0}$ in the complex structure $J$. Let $\{n_{1},\cdots,n_{20}\}$ be a basis of $Pic(S)$. Let $\alpha:=(\alpha_{1},\cdots,\alpha_{20})\in \mathrm{R}^{20}$. Define $\|\alpha\|:=\max_{1\leq i\leq 20} \{|\alpha_{i}|\}$. For any given positive number $\beta$  if $\|\alpha\|$ is small enough (depending on $\beta$) then the class \begin{equation}\omega_{J}^{\alpha,\beta}:= \beta\omega_{J}^{0}+\sum_{k=1}^{20}\alpha_{k}n_{k}\end{equation} is Kahler. First of all, $\omega_{J}^{\alpha,\beta}$ is real and of type $(1,1)$. By our assumption about $\beta$ and $\|\alpha\|$ we can assume that $(\omega_{J}^{\alpha,\beta})^{2}>0$. Therefore $\omega_{J}^{\alpha,\beta}$ is in a small enough neighborhood of $\beta\omega_{J}^{0}$ in the positive cone $V(S)^{+}$. Since the dimension of the positive cone is equal to the dimension of the Kahler cone and $\beta\omega_{J}^{0}$ is in the interior of the Kahler cone we know that a small enough neighborhood of  $\beta\omega_{J}^{0}$ in the positive cone is actually contained in the Kahler cone. So $\omega_{J}^{\alpha,\beta}$ is Kahler.

 $\omega_{J}^{\alpha,\beta}\in \langle p,q\rangle^{\perp}_{\mathrm{R}}:=(\langle p,q\rangle\otimes\mathrm{R})^{\perp}$ where $\langle p,q\rangle$ is  the sublattice generated by $p,q$ (which is the transcendental lattice of the singular $K3$ surface $S$). We also know that $pr(\omega_{J}^{\alpha,\beta})\in (U^{'}_{\mathrm{R}})^{\perp}$. Moreover $U^{'}_{\mathrm{R}}\subset\langle p,q\rangle^{\perp}_{\mathrm{R}}$. Therefore \begin{equation}pr(\omega_{J}^{\alpha,\beta})\in \langle p,q,f,\sigma_{0}\rangle^{\perp}_{\mathrm{R}}\end{equation}
We set $\omega_{J}=\omega_{J}^{\alpha,\beta}$. Then \begin{equation}\check{\omega}_{I}^{2}-\check{B}^{2}={1\over ((\omega_{J}^{\alpha,\beta})\cdot f)^{2}}(q^{2}+(Re(\tau))^{2}p^{2}-2Re(\tau)p\cdot q-(pr(\omega_{J}^{\alpha,\beta}))^{2})\end{equation}\begin{equation}D\cdot\check{B}={1\over (\omega_{J}^{\alpha,\beta})\cdot f}(pr(\omega_{J}^{\alpha,\beta})\cdot D)\end{equation}When $\delta$ belongs to type 1,2 or 3  we are going to show that the condition $Re(\exp(\check{B}+i\check{\omega}_{I}),\delta)=0$ is too strong to be satisfied by a generic perturbation of $\omega_{J}^{\alpha,\beta}$. For type 1, $s$ is a nonzero integer. If $D\cdot\check{B}$ happens to be a nonzero integer we just redefine $\omega_{J}^{\alpha,\beta}$ by adding $c\sigma_{0}$ (with a small enough positive constant $c$) to it. Then due to (98) the numerator of $D\cdot\check{B}$ does not change while the denominator is modified by adding a real number $c$. Hence for such a small perturbation of $\omega_{J}^{\alpha,\beta}$ the integral property $D\cdot\check{B}=a\ nonzero\ integer$ can not hold. Note that this means that $$Re(\exp(\check{B}+i\check{\omega}_{I}),\delta)\neq0$$ for all $\delta$'s of type 1. 

For type 2 or 3, $s,r, q^{2}+(Re(\tau))^{2}p^{2}-2Re(\tau)p\cdot q$ are fixed and rational. $r\check{\omega}_{I}^{2}-r\check{B}^{2}+2D\cdot\check{B}$ must be rational. We can assume that $pr(\omega_{J}^{\alpha,\beta})\neq0$. In fact if that is not the case we can perturb $\omega_{J}^{\alpha,\beta}$ along a direction in $(U^{'})^{\perp}\cap Pic(S)$. Then we will have $pr(\omega_{J}^{\alpha,\beta})\neq0$. Here $(U^{'})^{\perp}\cap Pic(S)\neq\emptyset$ because $\rho(S)$ is too large. Since $\langle p,q,f,\sigma_{0}\rangle$ has signature $(3,1)$ we know $\langle p,q,f,\sigma_{0}\rangle^{\perp}$ has signature $(0,18)$.   So $(pr(\omega_{J}^{\alpha,\beta}))^{2}\neq0$. Then we know that for a generic choice of $\omega_{J}^{\alpha,\beta}$ such that $(pr(\omega_{J}^{\alpha,\beta}))^{2}$ is irrational, $r\check{\omega}_{I}^{2}-r\check{B}^{2}$ and hence also $2D\cdot\check{B}$ must be nonzero. We also notice that the denominator of $D\cdot\check{B}$ is not a fixed multiple of the denominator of $\check{\omega}_{I}^{2}-\check{B}^{2}$. Therefore the same perturbation of Kahler class argument (adding $c\sigma_{0}$) shows that for some small perturbation of $\omega_{J}^{\alpha,\beta}$ the condition $Re(\exp(\check{B}+i\check{\omega}_{I}),\delta)=0$ is not true for all $\delta$'s of type 1, 2 or 3.

We still need to handle type 4 for which the condition $Re(\exp(\check{B}+i\check{\omega}_{I}),\delta)=0$ becomes \begin{equation}D\cdot\check{B}=0\end{equation}The above perturbation of Kahler class argument can fail in this case. It fails if and only if \begin{equation}D\cdot \eta=0 \ for\ any\ \eta\in\langle p,q,f,\sigma_{0}\rangle^{\perp}_{\mathrm{R}}\end{equation}Clearly if (102) is true then (101) always holds no matter how we perturb $\omega_{J}^{\alpha,\beta}$. And if (101) is true and there is an $\eta\in\langle p,q,f,\sigma_{0}\rangle^{\perp}_{\mathrm{R}}$ such that $D\cdot \eta\neq0$ then we can perturb $\omega_{J}^{\alpha,\beta}$ by adding $c\eta$ to get $D\cdot \check{B}\neq0$.

Now let us show that when $\delta$ is of type 4 the condition \begin{equation}(\exp(\check{B}+i\check{\omega}_{I}),\delta)=0\end{equation} and the condition (102) can not both hold. Clearly this  (together with the discussion on type 1,2 and 3) implies that for  some Kahler class $\omega_{J}^{\alpha,\beta}$ on $(S,J)$ we have \begin{equation}(\exp(\check{B}+i\check{\omega}_{I}),\delta)\neq0, \ for\ any\ \delta\in \Delta(\check{S})\end{equation} Then the theorem is proved.

Let us suppose that (103) and (102) both hold for some $\delta$ of  type 4. By computation  (103) is equivalent to \begin{equation}Re(\exp(\check{B}+i\check{\omega}_{I}),\delta)= D\cdot\check{B}=0\end{equation}\begin{equation} Im(\exp(\check{B}+i\check{\omega}_{I}),\delta)=D\cdot\check{\omega}_{I}=0\end{equation} The first equation is already implied by (102).  We also have the following obvious relations (since $D\in NS(\check{S})$)\begin{equation}Im(\check{\Omega}_{I})\cdot D=0 \ i.e.\ p\cdot D=0\end{equation}\begin{equation}Re(\check{\Omega}_{I})\cdot D=0\ i.e.\ ({1\over 2}(Im(\tau))^{2}p^{2}f+f+\sigma_{0})\cdot D=0\end{equation} Since $\langle p,q,f,\sigma_{0}\rangle^{\perp}$ has signature $(0,18)$, (102) implies that $D\in\langle p,q,f,\sigma_{0}\rangle$. Note that $U^{'}$ is the orthogonal complement of $\langle p, q\rangle$ in $\langle p,q,f,\sigma_{0}\rangle$. Then (107) and (106) together imply that $D=mf+n\sigma_{0}$ for some integers $m,n$. Then (108) becomes \begin{equation}m-n+{D_{p,q}\over 2p^{2}}n=0\end{equation}We have used $f^{2}=0$, $f\cdot\sigma_{0}=1$ and $\sigma_{0}^{2}=-2$. On the other hand we have $$D^{2}=2mn-2n^{2}=-2$$ which implies \begin{equation}n(m-n)=-1\end{equation} Since $m,n$ are integers we must have $n=\pm1$ and $m-n=\mp1$. Substitute these results back to (109) and we can get \begin{equation}D_{p,q}=2p^{2}\end{equation} But this violates our assumption. This is a contradiction.\\

Finally we want to point out that the condition that $\sqrt{D_{p,q}}$ is integral and the condition $D_{p,q}\neq2p^{2}$ can be satisfied at the same time. The following positive definite even integral quadratic form $$Q_{p,q}:=\left(\begin{array}{cc}2 & 0\\ 0 & 8\end{array}\right)$$is such an example. $\diamondsuit$\\

 Theorem 3.3 and theorem 6.3 tell us that there is a covering map \begin{equation}\pi^{-1}(\mathcal{P}_{0}(\check{S}))\rightarrow \mathcal{P}_{0}(\check{S})\end{equation} such that there are stability conditions in $\pi^{-1}(\mathcal{P}_{0}(\check{S}))$ whose central charge  is \begin{equation}Z(E)=(\exp(\check{B}+i\check{\omega}_{I}),v(E))\end{equation}Let $Stab^{*}(\check{S})$ be the good component containing these stability conditions. We will write $Z(E)$ as $Z(l), l=v(E)$ to emphasize that the central charge depends only on the Mukai vector $l\in \mathcal{N}(\check{S})$ and not on $E$ directly. Clearly $Z(l)$ can be extended by the same formula (113) to be defined for all  $l\in H^{*}(\check{S},\mathrm{Z})$.

 \begin{dfn} Let $v_{i},v_{j}\in H^{*}(\check{S},\mathrm{Z})$  and $\tilde{\gamma}=(i,j)$. We define the $generalized\ stability\ wall$ $\tilde{W}_{\tilde{\gamma}}$\begin{equation}\tilde{W}_{\tilde{\gamma}}=\{\sigma=(Z,\mathcal{P})\in Stab^{*}(\check{S}):Z(v_{i})/Z(v_{j})\in\mathrm{R}_{>0}\}\end{equation}\end{dfn}Although the definition of $\tilde{W}_{\tilde{\gamma}}$ is similar to the definition of the stability wall $W_{\tilde{\gamma}}$ in section 3 there is a difference. In the definition of $\tilde{W}_{\tilde{\gamma}}$ we say nothing about  objects in $\mathcal{D}(\check{S})$ while in the definition of $W_{\tilde{\gamma}}$ we require that the elements $v_{i},v_{j}\in \mathcal{N}(\check{S})$ are Mukai vectors of some objects in $\mathcal{D}(\check{S})$. Every stability wall is a generalized  stability wall.  Of course if $v_{i},v_{j}\in\mathcal{N}(\check{S})$ are Mukai vectors of some objects used in the definition 3.2 then $\tilde{W}_{\tilde{\gamma}}$ is a stability wall. \\

Let $\{l_{1},\cdots, l_{20}\}$ be a basis of $Pic(S)=Pic(S,J)$. By definition we have $\{\mu(l_{i})\in H^{*}(\check{S},\mathrm{Z})\},1\leq i\leq20$. Let $\tilde{\Gamma}$ be the set of pairs $(i,j), 1\leq i,j\leq20$ with $i<j$. Then for any $\tilde{\gamma}=(i,j)\in \tilde{\Gamma}$  we have the generalized stability wall\begin{equation}\tilde{W}_{\tilde{\gamma}}=\{\sigma=(Z,\mathcal{P})\in Stab^{*}(\check{S}):Z(\mu(l_{i}))/Z(\mu(l_{j}))\in\mathrm{R}_{>0}\}\end{equation}

\begin{thm} We make the assumptions of theorem 6.3. Then there is some Kahler class $\omega_{J}$ on $(S,J)$ such that (89) is true and  $Re Z(\mu(l_{i}))\neq0$ for any $i$. We  change $l_{i}$ to $-l_{i}$ for some $i$'s if necessary (see the proof for its meaning). Then the projections of  generalized stability walls $\{\pi(\tilde{W}_{\tilde{\gamma}})\}, \tilde{\gamma}\in\tilde{\Gamma}$ from the space of  stability conditions on $\check{S}$ to $\mathcal{N}(\check{S})\otimes\mathrm{C}$ intersect at $\check{B}+ i\check{\omega}_{I}$.  \end{thm}

\noindent {\bf Proof}  First we show that $Re Z(\mu(l_{i}))\neq0$ for any $i$.

 We already know $Z(\mu(f))=1$. For $Re Z(\mu(\sigma_{0}))$ we see that $\mu(\sigma_{0})=(1,0,1)$ is of type 2 in the proof of theorem 6.3. By the proof of theorem 6.3 there is an  $\omega_{J}$ such that both (89) and $Re Z(\mu(\sigma_{0}))\neq0$ are true. For $l_{i}\in\Gamma^{'}$ (also note that $l_{i}\in Pic(S)$), we have $Re Z(\mu(l_{i}))=l_{i}\cdot\check{B}$ which is the same as (96) with $D$ replaced by $l_{i}$. But here we know that $l_{i}\in \langle p,q,f,\sigma_{0}\rangle^{\perp}$. So the condition (102) is not true. But then the proof of theorem 6.3 tells us that $Re Z(\mu(l_{i}))\neq0$ for some small perturbation of $\omega_{J}$. In general $l_{i}=al_{i}^{a}+mf+n\sigma_{0}$ where $l_{i}^{a}\in\Gamma^{'}$ and $a,m,n$ are integers. Again a perturbation of Kahler class guarantees that $Re Z(\mu(l_{i}))\neq0$ for all $i$ and (89) holds. \footnote{As in the proof of theorem 6.3 the case of $l_{i}\cdot \check{B}\neq 0$ requires a different treatment from other cases. So we have discussed it before the general cases.}

 By theorem 6.2 every  $Z(\mu(l_{i}))$ is real i.e. $Z(\mu(l_{i}))=Re Z(\mu(l_{i}))$. So after changing some $l_{i}$ to $-l_{i}$ for some $i$'s we can assume that $\arg Z(\mu(l_{i}))>0$ for any $i$. We know that $\exp(\check{B}+ i\check{\omega}_{I})\in \mathcal{P}_{0}(\check{S})$ and therefore is in $\pi(Stab(\check{S}))$. Clearly it is in $\pi(Stab^{*}(\check{S}))$ for some connected component $Stab^{*}(\check{S})$. So we can define generalized stability walls $\{\tilde{W}_{\tilde{\gamma}}\}, \tilde{\gamma}\in\tilde{\Gamma}$. The rest of the theorem is only a restatement of the fact: $\arg Z(\mu(l_{i}))= \arg Z(\mu(l_{j}))$ for each pair $\tilde{\gamma}=(i,j)$. $\diamondsuit$\\

\begin{note}Theorem 6.4 and theorem 5.1 together is a statement about a correspondence of mirror symmetry in an attractor background. On both sides we have correspondent conditions about  central charges (i.e. alignment of many central charges with charges related by the mirror map). On one side of this correspondence we  get intersections of generalized stability walls in the space of stability conditions.   The only remaining part to get a complete  correspondence  of type (5) is defining stability conditions on derived Fukaya categories properly.\end{note}
\begin{note}We may wonder if one can find a statement in correspondences of type  (3). This problem is harder. It  requires us to verify that stable objects are related by the mirror symmetry (this is much more refined than homological mirror conjecture) but mirror symmetry of $K3$ surfaces used in this paper is mapping cohomology classes to cohomology classes (this is coarser than homological mirror conjecture). The following question should be answered for the study of correspondences of type (3):

When is a  generalized stability wall in theorem 6.4 a stability wall? \end{note}
\begin{note}It is natural to ask if there is a generalization of main theorems to $S\times T^{2}$. After all this is where the story begins. The only reason that we do not study this problem  is that people do not know how to build stability conditions on Calabi-Yau threefolds (but see \citep{BMT}). \end{note}
\begin{note}In the paper we have directly verified that $Z_{K3}(l)$ equals $Z(\mu(l))$ up to a real normalization when the complex moduli for the former and the complexified Kahler moduli for the latter are related by mirror symmetry (of course the former is also assumed to be an attractor background). This kind of statements have been widely used by physicists (for more general settings) and can lead to nontrivial mathematical predictions (see \citep{DRY}). However the author feels that  it is not clear why this should be true from a purely mathematical perspective. It should be considered as a prediction of mirror symmetry and requires justifications.\end{note}

\begin{note} One can speculate if a mirror $\check{S}$  of a $K3$ surface in an attractor background is also in some sort of attractor background. This would make mirror symmetry more perfect. Of course for $\check{S}$ the moduli must be complexfied Kahler moduli. In fact there does exist an attractor mechanism in the space of complexifed Kahler moduli for black holes in  another type of string theory. A lot of papers have been written about this attractor mechanism. However the central charges in this attractor mechanism receive corrections that are not under control when we are far from the large volume limits. So even though it might be true that the mirror of an attractor background (in the complex moduli) is an attractor background (in the complexified Kahler moduli) it is not clear how to verify it.\end{note}
\begin{note}The wall crossings we study in this paper are static in the sense that we are sitting on attractor backgrounds where many stability walls intersect. A more dynamical approach would be starting from somewhere away from attractor backgrounds and then move to an attractor background. During this process we might encounter some stability walls. If we move along the attractor flow (which seems to be the most natural choice) then we can actually derive some very sophisticated wall crossing formulas. One such example is \citep{DeM}. This point of view may be very helpful when we try to formulate and prove correspondences of  (3).\end{note}
\begin{note} Readers may have realized that black holes actually do not play a very essential role in this paper. All we need from them is the motivation to study singular $K3$ surfaces and the idea that we should compute the central charge. However if one wants to use the dynamical approach in the previous remark to study wall crossings of stable objects there is a way to make the black hole perspective essential. Instead of just playing in the space of complex moduli (or complexified Kahler moduli) one can study the existence of black hole solutions and use it to deduce nontrivial wall crossing information. For a sample of this big subject see \citep{DeM}. This is very puzzling for a geometer: how can the existence and counting of special Lagrangians and stable objects in $\mathcal{D}(X)$ be related to black holes? Clearly we have just seen a tip of a huge iceberg. We hope  to discuss this topic in the future.\end{note}

\end{document}